\documentclass[12pt]{article}
\usepackage{amsmath,amsfonts,amssymb,amsthm}
\usepackage{enumerate}
\usepackage{thm-restate}
\usepackage{fullpage}
\usepackage{cite}

\usepackage[margin=1in]{geometry}

\theoremstyle{plain}
\newtheorem{theorem}{Theorem}[section]

\newtheorem{proposition}[theorem]{Proposition}
\newtheorem{corollary}[theorem]{Corollary}
\newtheorem{lemma}[theorem]{Lemma}
\theoremstyle{definition}
\newtheorem{definition}[theorem]{Definition}
\newtheorem*{remark}{Remark}
\newtheorem*{observation}{Observation}

\newcommand{\N}{\mathbb N}
\newcommand{\Z}{\mathbb Z}
\newcommand{\Q}{\mathbb Q}
\newcommand{\Gal}{\operatorname{Gal}}
\newcommand{\Col}{\operatorname{Col}}
\newcommand{\sgn}{\textrm{sgn}}
\newcommand{\co}{\textrm{co}}
\newcommand{\NP}{\mathcal{N}}
\newcommand{\OK}{\mathcal O_K}
\newcommand{\pp}{\mathfrak p}

\title{A Dynamical Analogue of Sen's Theorem}

\author{Mark O.-S. Sing}

\begin{document}
\maketitle

\begin{abstract}
	We study the higher ramification structure of dynamical branch extensions, and propose a connection between the natural dynamical filtration and the filtration arising from the higher ramification groups: each member of the former should, after a linear change of index, coincide with a member of the latter. This is an analogue of Sen's theorem on ramification in $p$-adic Lie extensions. By explicitly calculating the Hasse-Herbrand functions of such branch extensions, we are able to show that this description is accurate for some families of polynomials, in particular post-critically bounded polynomials of $p$-power degree. We apply our results to give a  partial answer to a question of Berger~\cite{BergerIterated} and a partial answer to a question about wild ramification in arboreal extensions of number fields~\cite{AitkenHajirMaire,BIJJLMRS}.
\end{abstract}

\section{Introduction}
Many guiding questions in arithmetic dynamics arise from or are inspired by analogies to well-studied objects in arithmetic geometry. Here, we formulate a tentative dynamical analogue of Sen's theorem, and prove it in certain cases. Sen noticed that, for Galois extensions whose Galois groups are $p$-adic Lie groups, there is a remarkable connection between the $p$-adic Lie filtration, which depends only on the Lie group, and the filtration by upper ramification subgroups: the two mutually refine each other in a precise way after a linear change of index~\cite{Sen}. In our dynamical setting, we replace $p$-adic Lie groups and the Lie filtration with  ``branch extensions'' and their ``branch filtration'' (see Section 1.1 for definitions and notation). For those familiar with arboreal representations, we are taking a single branch of the tree, filtered by height up the branch.

Our dynamical version of Sen's theorem says that, after possibly extending the ground field and making a linear change in index, each member of the branch filtration coincides exactly with a member of the upper ramification filtration. The upper ramification filtration is in general quite difficult to understand, and captures subtle arithmetic phenomena, while the branch filtration is quite simple and dynamically natural: starting from our ground field $K$, we have a tower of extensions $K_n$ over $K$ obtained by adjoining a compatible sequence (``branch'') of preimages of the base point. We are able to give a general sufficient criterion for our result to hold: it applies to extensions associated to so-called ``tamely ramification-stable'' branches. In our situation, ``tamely'' simply means that $p$ does not divide a certain quantity $d$, which is the limiting valuation of the members of the branch. Such branches are particularly striking from a dynamical perspective, exhibiting a kind of stability in the structure of their higher ramification: the intermediate Hasse-Herbrand functions associated to $K_n/K_{n-1}$ are identical up to according to small and well-controlled errors, neglecting scaling. For these branches, we obtain our main result:
\begin{restatable*}{theorem}{TheoremMain}\label{TheoremMain}
	Suppose our branch, associated to the polynomial $P(x)$ and base point $\alpha_0$, is tamely ramification-stable over $K$. Then $K_\infty/K$ is arithmetically profinite, and there is a constant $V$ such that for all $n$,
	$$K_n = K_\infty^{((V-1)n + 1)}.$$
\end{restatable*}

A more literal, and weaker, restatement of Sen's theorem in the dynamical setting would be that the two filtrations refine each other, again, after a linear change of index. However, for one of our applications, to a question of Berger~\cite{BergerIterated}, we need this stronger formulation.

We are able to give a general sufficient criterion for a branch to be tamely ramification-stable, Proposition~\ref{PropComposition}. This criterion consists of two pieces: that $p$ does not divide $d$, and verifying an inequality depending only on the valuations of the coefficients of $P(x)$ and the valuation of $\alpha_0$. Some branches which are not tamely ramification-stable may become so after extending the ground field and re-indexing the branch; we call such branches \textit{potentially} tamely ramification-stable.

Using this criterion, we are able to show that if $P(x)$ is either post-critically bounded or prime degree, and we take a branch such that $p$ does not divide the associated constant $d$, then it is \textit{potentially} tamely ramification-stable, and use this information to characterize higher ramification in the associated extension:

\begin{restatable*}{corrollary}{CorollaryMain}\label{CorollaryMain}
	Let $P(x)$ be a polynomial which either has degree $p$, or is post-critically bounded and has degree $p^r$. Take any nontrivial branch for $P(x)$, and suppose $p$ does not divide the constant $d$ associated to the branch.
	
	Then the dynamical branch extension $K_\infty/K$ is arithmetically profinite, and there are constants $N$ and $V$ such that after replacing $K$ by $K_N$, 
	$$K_n = K_\infty^{((V-1)(n-N) + 1)},$$
	for all $n$.
\end{restatable*}

For any particular branch, it is not difficult to apply our criteria to check whether or not it is (potentially) tamely ramification-stable, so long as one knows that $p$ does not divide $d$. In fact, our criterion is almost entirely effective: only the stipulation that $p$ does not divide $d$ is not known to be effective. Each branch determines certain ``limiting ramification data'' from which one can completely recover the Hasse-Herbrand function of the associated branch extension in the tamely ramification-stable case when $d$ is known. The calculation of the limiting ramification data depends only on $P(x)$ and some of the initial entries of the branch (the number of entries needed is itself effective). While we lack a general algorithm to determine $d$, it can be calculated in many particular instances.

We apply our results to provide a partial to answer two questions. One is raised by Berger~\cite{BergerIterated}, who asks: is it possible to show by elementary methods that if $K_\infty/K$ is Galois and the base point is a uniformizer then its Galois group is abelian? This is known to be true by Berger~\cite{BergerLifting} using quite sophisticated methods from $p$-adic Hodge theory. Our main theorem involves more elementary tools, and proves allows us to re-prove this fact in some situations:
\begin{restatable*}{theorem}{ThmBergerApp}\label{ThmBergerApp}
	Assume $p$ is odd. Suppose $\alpha_0$ is a uniformizer for $K$, $P'(0)$ is nonzero, and we are given a branch associated to $P(x)$ and $\alpha_0$ which is tamely ramification-stable.
	
	If $K_\infty/K$ is Galois, it is also abelian.
\end{restatable*}

The other question is suggested by both Aitken, Hajir, and Maire (Question 7.1 in~\cite{AitkenHajirMaire}) and Bridy, Ingram, Jones, Juul, Levy, Manes, Rubinstein-Salzado, and Silverman (Conjecture 6 in~\cite{BIJJLMRS}), who essentially ask if it is possible for an arboreal extension over a number field to be ramified at finitely many primes \textit{but not wildly ramified}. It turns out that this is not possible for polynomials of prime-power degree:
\begin{restatable*}{theorem}{BIJJLMRSApplication}\label{BIJJLMRSApplication}
	Let $F$ be a number field and $\pp$ a prime of $F$ lying over a rational prime $p$. Let $P(x) \in \mathcal O_F[x]$ be a monic polynomial of degree $p^r$ such that $P(x) \equiv x^{p^r} \mod \pp$, and let $\alpha_0 \in F$.
	
	Then the arboreal representation associated to $P(x)$ and $\alpha_0$ is infinitely wildly ramified.
	
	If, further, $P(x)$ has prime degree and $v(\alpha_0) \neq 0$, or is post-critically bounded with no restriction on $v_\pp(\alpha_0)$, and there is a branch over $\alpha_0$ whose associated constant $d$ is not divisible by $p$, then every higher ramification subgroup over $\pp$ of the arboreal representation is nontrivial.
\end{restatable*}

\subsection{Prior work.}
While we have stated our main result as an analogy to Sen's theorem, our initial motivation comes from arithmetic dynamics and the structure of arboreal representations associated to post-critically finite maps. Arboreal representations, first introduced by Odoni~\cite{Odoni}, have been a subject of significant focus in arithmetic dynamics. This recently culminated in the resolution of (one version of) Odoni's conjecture over number fields, in prime degree by Looper~\cite{Looper}, in all even degrees and certain odd degrees by Benedetto and Juul~\cite{BenedettoJuul}, and finally for all degrees by Specter~\cite{Specter}. The branch extensions we tackle are the subextensions of the full arboreal extension which are associated to a single branch of the full preimage tree. The extensions we study appear within the full arboreal representation and the ramification along such branches is quite important to the aforementioned results on Odoni's conjecture. Additionally, Andrews and Petsche~\cite{PetscheAndrews} as well as Ferraguti and Pagano~\cite{FerragutiPagano} have also used ramification information to prove interesting results about \textit{abelian} arboreal representations over number fields. Our results are finer than necessary for any of the papers mentioned, but the important role ramification plays in those results suggests the potential value of the more detailed and delicate ramification information that we obtain. Though arboreal extensions over global fields are still quite mysterious, even less is known over local fields. Recently Anderson, Hamblen, Poonen, and Walton~\cite{AHPW} studied full arboreal extensions in the local setting for polynomials of the form $x^n + c$. In fact, they even produce an example which shows that a literal dynamical analogue of Sen's theorem cannot hold in full generality, even in the case of prime degree.

The case of post-critically bounded polynomials is of particular dynamical interest because it includes the post-critically finite polynomials of prime-power degree. Currently, the arboreal representations of post-critically finite polynomials are not well-understood, but it is known that they have arboreal representations which are ramified at only finitely many primes~\cite{AitkenHajirMaire}, so one would expect their arboreal representations to largely be controlled by their local behavior at those primes. Our result reveals initially unexpected structure to their wild ramification at the prime in question.

Some other work has been done with extensions of the kind we consider. Both Berger~\cite{BergerIterated} and Cais and Davis~\cite{CD} study them (under the name ``$\phi$-iterate extensions'') with the machinery of $p$-adic Hodge theory, and show that if these extensions are Galois they must be abelian. Cais, Davis, and Lubin~\cite{CDL} study the ramification in a somewhat more general setting, using similar methods to ours to give a characterization of arithmetically profinite extensions -- it is an important corollary of Sen's theorem that $p$-adic Lie extensions are arithmetically profinite. The dynamical case of their result applies to a broader class of polynomials than ours, with the restriction that the base point is a uniformizer. For the polynomials considered in this paper, we are able to relax this restriction on the base point and obtain more precise information about the ramification of our extensions.

\subsection{Outline of the paper.}

The structure of our paper is as follows. Section 2 establishes some preliminary facts, including weaker descriptions of ramification in branch extensions. Section 3 uses the preliminaries of Section 2 to furnish more detailed ramification information, then introduces the notion of (potentially) tamely ramification-stable extensions, ending with a proof of our main result. Section 4 discusses the two aforementioned applications to the questions of Berger and Bridy, Ingram, Jones, Juul, Levy, Manes, Rubinstein-Salzado, and Silverman, and also the effectivity of our results for which we provide algorithms that can sometimes allow to us to verify that a given dynamical system satisfies the assumptions of our theorems. We apply these algorithms to provide an example of a tamely ramification-stable branch which is not prime degree or post-critically bounded.

\subsection{Notation.}
\begin{enumerate}[\indent--]
	\item $p$ is a prime,
	\item $K$ is a valued field of characteristic zero complete for a discrete valuation and with perfect residue field of characteristic $p$. For example, a finite extension of $\Q_p$ or of $\widehat{\Q_p^{ur}}$,
	\item $\bar K$ is a fixed algebraic closure of $K$,
	\item $\OK$ is the ring of integers of $K$ and $\pi_K$ a uniformizer of $\OK$,
	\item $P(x) \in \OK[x]$ is a monic polynomial of degree $q=p^r$ such that $P(0) = 0$ and $P(x) \equiv x^q \mod \pi_K$,
	\item $v$ is a valuation for which $\bar K$ is complete, such that $v(p)$ and the valuations of the coefficients of $P$ are integers, and there is a subfield $E$ of $K$ such that $[K:E]$ is finite and $v(E) = \Z$,
	\item $P^n(x)$ denotes the $n$th iterate of $P$,
	\item $\alpha_0 \in K$ is the base point, and we assume $v(\alpha_0) \neq 0$,
	\item $e_{K/E}$ is the ramification index of the extension $K/E$,
	\item $(\alpha_n)_{n\in \N}$ is a sequence in $\bar K$ such that $P(\alpha_n) = \alpha_{n-1}$ and not all entries are zero (such sequences may also be called \textit{branches}, in keeping with the arboreal nomenclature),
	\item $K_n=K(\alpha_n)$, $K_0 = K(\alpha_0) = K$ and $K_\infty=\bigcup_n K_n$,
	\item $\pi_n$ are uniformizers of $\mathcal O_{K_n}$, along with units $u_n\in \mathcal O_{K_n}$, and integers $d_n$ such that $\alpha_n = u_n\pi_n^{d_n}$,
	\item $d=\lim_{n\rightarrow\infty} d_n$, if this limit exists,
	\item $\NP_n$ is the Newton polygon of $P(x+\alpha_n) - \alpha_{n-1}$, and $\co\NP_n$ is the copolygon,
	\item $\phi_n$ is the Hasse-Herbrand function for $K_n/K_{n-1}$,
	\item $\Phi_n$ is the Hasse-Herbrand function for $K_n/K$,
	\item $\Gamma_K$ is the absolute Galois group of $\bar K$ over $K$,
	\item $\Gamma_K^\nu$ the subgroup associated to a nonnegative real number $\nu$ by the upper-numbering of the higher ramification subgroups,
	\item $b_m$ the $m$th ramification break (see below).
\end{enumerate}

One may take $E=K$ with $v$ an appropriately normalized valuation, but we separate $E$ from $K$ here in anticipation of changing the base field $K$. The choice of $E$ amounts to fixing a choice of valuation on $\bar E = \bar K$ which doesn't vary as we vary $K$.

At times we treat the cases $v(\alpha_0) > 0$ and $v(\alpha_0) < 0$ separately, and refer to them as the \textit{integral} and \textit{non-integral} cases, respectively. It is easy to see from the Newton polygon that $v(\alpha_n)$ has the same sign as $v(\alpha_0)$ for all $n$.

By \textit{conjugation} we mean conjugation by a nonconstant linear map, rather than the more typical (for dynamics) conjugation by a linear fractional transformation, because we work with polynomials.

For a polynomial $f(x)$, we denote by $f_i$ the coefficient of $x^i$.

We assume that the reader is familiar with local fields and higher ramification. Serre~\cite{Serre} covers much of this (Chapter IV) but we do \textit{not} assume that our extensions are Galois, which at times requires slightly different tools; fortunately, Lubin has collected these in an excellent expository article~\cite{Lubin}. One should take caution in passing between these sources: their ramification numberings differ, and in the present paper we adopt Lubin's numbering. The reader should have a least a passing familiarity with the notion of arithmetically profinite extensions, such as in Wintenberger~\cite{Wintenberger}.

When $L/K$ is arithmetically profinite, we denote its $m$th ramification break by $b_m$ and the $m$th elementary subfield of $L$ by $L^{(m)}$. This is the subfield of $L$ which is fixed by $\Gamma_K^{b_m}$. By convention, we set $K_\infty^{(m)} = K$ when $m$ is nonpositive.

We remark that many of the results of this paper hold in equicharacteristic $p$ in the presence of separability assumptions. In many cases the proofs become simpler in this situation. However, an important fact we prove about post-critically bounded polynomials (that they have potentially good reduction) is no longer true.

\section{Preliminaries}

The results here are used primarily as tools for our main theorems. However, some may be of independent interest, which we have tried to mark as propositions.

Any post-critically bounded (PCB) polynomial of $p$-power degree can, after possibly extending the ground field, be put in the same form as the polynomials we consider: monic, with integral coefficients, and fixing zero. In fact, after conjugation, a post-critically bounded polynomial satisfies even stronger constraints on its coefficients. A few other versions and proofs of Proposition~\ref{PropPCB} appear in the literature~\cite{Anderson,Epstein,BJL}.
\begin{proposition}\label{PropPCB}
	If a polynomial has degree $q=p^r$ and is post-critically bounded, then it has a conjugate $P(x)$ which is monic, integral, fixes $0$ and further satisfies
	$$v(P_i) + v(i) \geq v(q) = rv(p)\ \ \ \textrm{ for all } 1\leq i \leq q.$$
\end{proposition}
\begin{proof}
	Let $R(x)$ be the given polynomial. After conjugating, we may assume that $R(x)$ is monic and fixes zero; conjugates also remain post-critically bounded. This conjugation may require taking a $(p-1)$th root of the leading coefficient of $R(x)$ and adjoining a fixed point of $R$ to the ground field. Call this conjugate $P(x)$. It does not necessarily have integral coefficients at this point, but we will show that $\frac{P'(x)}{q}$ is in $\OK[x]$, from which the claimed inequality follows, and hence integrality as well.
	
	Suppose otherwise, that $\frac{P'(x)}{q}$ is not in $\OK[x]$. This guarantees a positive slope in the Newton polygon for $\frac{P'(x)}{q}$, the steepest slope of which ends at the vertex associated to the leading term. This slope must be strictly steeper than the steepest slope of the Newton polygon of $P(x)$ because every non-leading vertex moves down in passing from $P(x)$ to $\frac{P(x)}{q}$. However, this means if we take a critical point associated to this steepest slope, $v(P(c)) = qv(c) < v(c)$, hence $v(P^2(c)) < q^2v(c)$ and so on, so $v(P^n(c)) \rightarrow -\infty$ and hence the critical orbit is unbounded.
\end{proof}

This only tells us that a post-critically bounded polynomial has some conjugate of the desired form. Conjugation moves the base point, and a priori could leave us with a base point of valuation zero, contrary to our requirements. It turns out that, after possibly replacing $P$ by an iterate, there is always a choice of conjugate such that the new base point has nonzero valuation. This is elaborated on in Section 4.2.

Eventually, we will look at $\NP_n$ and $\co\NP_n$, the Newton polygon and copolygon associated to $P(x+\alpha_n) - \alpha_{n-1}$. When we expand this expression, the coefficients of the resulting polynomial involve binomial coefficients, and so to control these polygons we need some control over the binomial coefficients as well.
\begin{lemma}\label{LemmaBinomialDrop}
	Fix a positive integers $i,j,k$ with $j \geq i$ and $j\geq p^k$.
	\begin{enumerate}[(i)]
	\item If $p^k \leq i < p^{k+1}$, then
	$$v\binom j{i} \geq v\binom j{p^k}.$$
	\item Additionally,
	$$v \binom j {p^{k+1}} \geq v\binom j {p^k} - v(p),$$
	with equality if and only if $v\binom j{p^k} \neq 0$.
	\end{enumerate}
\end{lemma}
\begin{proof}
	Both claims follow from Kummer's theorem~\cite{Kummer}, which states that the $p$-adic valuation of a binomial coefficient $\binom j i$ is $cv(p)$, where $c$ is the number of carries when adding $i$ and $j-i$ in base $p$. 
	
	Applying that theorem, we see that a lower bound for the valuation of $\binom j i$ when the leading base $p$ digit of $i$ is in the $\ell$th place is the number of consecutive zeros in the base-$p$ expansion of $j$ starting at the $\ell$th digit. Notice that if $i=p^\ell$ then this is exact, but it can be larger in general, from carries that occur	 before the $\ell$th digit.
	
	The condition $p^k\leq i < p^{k+1}$ says exactly that $i$'s leading base $p$ coefficient is in the $k$th place.
	
	From these observations, $(i)$ and the inequality of $(ii)$ are immediate by taking $\ell=k$ and $\ell=k+1$. As to the last claim: the quantities in question are integers, so equality is impossible if $v\binom j {p^k}$ is zero, and conversely if $v\binom j {p^k}$ is nonzero then the change from $p^k$ to $p^{k+1}$ loses exactly one of the aforementioned zeros in its base $p$ expansion.
\end{proof}

The next proposition is our first dynamical result, a kind of ramification semi-stability, although much less refined than our main result.
\begin{proposition}\label{PropSemistable}
	Let $P(x)$ and $(\alpha_n)$, be a branch and polynomial over a field $K$, as described in Section 1.1 (in particular, $P(x) \equiv x^q \mod \pi_K$ and the branch is nontrivial and based at an element of nonzero valuation). Then for all sufficiently large $n$:
	\begin{enumerate}[\indent (a)]
	\item $v(\alpha_{n+k}) = \dfrac{v(\alpha_n)}{q^k}$,
	\item the sequence $(d_{n+k})_{k\in\N}$ is constant, hence $\lim_{n\rightarrow \infty} d_n = d$ exists,
	\item $K_n/K_{n-1}$ is totally ramified of degree $q$.
\end{enumerate}
\end{proposition}
\begin{proof}
	Consider the polynomial $P(x) - \alpha_{n-1}$, of which $\alpha_n$ is a root. We may, by taking $n$ large enough, assume $\alpha_{n-1}$ is nonzero. Inspecting the Newton polygon of $P(x) - \alpha_{n-1}$, we see that, in the integral case where $v(\alpha_n)> 0$,
	
	$$v(\alpha_n) \leq \max \{v(\alpha_{n-1}) - 1, v(\alpha_{n-1})/2\}.$$
	
	Thus, for $n$ large enough, we obtain $v(\alpha_n) < 1$, in which case the Newton polygon of $P(x)-\alpha_n$ has just one segment, whose slope is $\frac{v(\alpha_n)}{q}$, which is also less than $1$. Hence $v(\alpha_{n+1}) = v(\alpha_n)/q$, and inductively this yields $v(\alpha_{n+k}) = v(\alpha_n)/q^k$.
	
	In the non-integral case, the Newton polygon immediately has a single slope, which gives us
	$$v(\alpha_{n+1}) = \frac{v(\alpha_{n})}{q} <0,$$
	and again the claim follows inductively. Together, these two cases verify $(a)$.
	
	Next, let $e_n$ be the ramification index of $K_n/K_{n-1}$. Take $n-1$ large enough that $(a)$ holds, and so we have both
	$$v(\alpha_n) = v(u_n\pi_n^{d_n}) = d_nv(\pi_n) = \frac{d_n v(\pi_{n-1})}{e_n},$$
	and
	$$v(\alpha_n) = \frac{v(\alpha_{n-1})}{q} = \frac{v(u_{n-1}\pi_{n-1}^{d_{n-1}})}{q} = \frac{d_{n-1}v(\pi_{n-1})}{q}.$$
	
	Comparing the two yields the following relation:
	\begin{equation}\label{dnenrelation}
	d_n= \frac{e_n}{q}d_{n-1}.
\end{equation}
	
	From (\ref{dnenrelation}), we see that if $e_n = q$, then $d_n = d_{n-1}$, and so to prove both $(b)$ and $(c)$ it suffices to verify $e_n=q$ for $n$ large enough. Evidently $e_n\leq q$, so we wish to show that this inequality is strict at most finitely often. Indeed, each time the inequality is strict, the $p$-adic valuation of $d_n$ is strictly smaller than the $p$-adic valuation of $d_{n-1}$. Since the sequence of valuations $(v_p(d_n))$ is a sequence of nonnegative integers, these strict drops can happen only finitely many times, and hence it is eventually constant.
\end{proof}

Using Proposition~\ref{PropSemistable}, we are able to take a step towards more precise information about the Newton polygons $\NP_n$.
\begin{lemma}\label{LemmaNP}
	For $n$ sufficiently large, the Newton polygon $\NP_n$ of $P(x+\alpha_n) - \alpha_{n-1}$ has at most $r+1$ vertices, whose $x$-coordinates can only be powers of $p$.
	
	Thus $\NP_n$ is the lower convex hull of the points $(p^k,y_{p^k})$, where the height $y_{p^k}$ is given by
	$$y_{p^k} = \min_{p^k\leq j \leq q} \left\{v\binom j {p^k} + v(P_j) + (j-p^k)v(\alpha_n)\right\}.$$
\end{lemma}
\begin{proof}
	Let $Q(x) = P(x+\alpha_n) - \alpha_{n-1}$. Expanding and collecting terms, we see that
	$$Q_i = \sum_{j=i}^q \binom j i P_j \alpha_n^{j-i}.$$
	Hence
	\begin{equation}\label{npheighteqn}
		v(Q_i) \geq \min_{i\leq j\leq q} \left\{v\binom j i + v(P_j) + (j-i)v(\alpha_n)\right\}
	\end{equation}
	
	The fractional parts of the terms in the minimum, which come from $(j-p^k)v(\alpha_n)$, are all distinct so long as $0<|v(\alpha_n)| \leq \frac 1 q$, and from Proposition~\ref{PropSemistable} we know this is the case for all sufficiently large $n$. As such, the terms themselves are distinct and so the inequality (\ref{npheighteqn}) is actually an equality. Additionally, $v(Q_1)\neq \infty$ since the minimum above is evidently finite.
	
	Since $Q_0 = 0$, but $Q_1 \neq 0$, the Newton polygon has a vertical line through $(1,v(Q_1))$. The leading coefficient is $1$, so there is another vertex at $(q,0)$.
	
	To show that $\NP_n$ only has vertices at prime powers, we will prove something slightly stronger: that $v(Q_i)$ for $i$ between $p^k$ and $p^{k+1}$ has valuation at least $v(Q_{p^k}) + (p^k - i)v(\alpha_n)$, or, in other words, such points $(i,Q_i)$ are above the line through $(p^k,v(Q_{p^k}))$ with slope $-v(\alpha_n)$. Because $|v(\alpha_n)|\leq \frac 1 q$, the slope of that line through $(p^k,v(Q_{p^k}))$ is so shallow, that this line always passes above $(q,0)$ and so no point above this line can be a vertex \textit{whether or not $(p^k,Q_{p^k})$ is itself a vertex}. Since we will prove that every point strictly between $p^k$ and $p^{k+1}$ does lie above such a line, none of them can be vertices, hence the only admissible locations for vertices are at prime powers.
	
	And so we compute, for $p^k \leq i <p^{k+1}$:
	\begin{align}
	v(Q_i) &= \min_{i\leq j \leq q}  \left\{v\binom j i + v(P_j) + (j-i)v(\alpha_n)\right\}\notag\\
	    &= \min_{i\leq j \leq q}  \left\{v\binom j i + v(P_j) + (j-p^k)v(\alpha_n)\right\} + (p^k - i)v(\alpha_n)\notag\\
	    &\geq \min_{p^k\leq j \leq q}  \left\{v\binom j i + v(P_j) + (j-p^k)v(\alpha_n)\right\} + (p^k - i)v(\alpha_n)\label{ineqbreakslope}\\
	\intertext{This is nearly the desired inequality, but with $v\binom j i$ rather than $v\binom j {p^k}$. To resolve this issue, we apply Lemma~\ref{LemmaBinomialDrop}, which tells us that if $p^k \leq i < p^{k+1}$, then
	$$v \binom j i \geq v \binom j {p^k}.$$
	Continuing where we left off at (\ref{ineqbreakslope}):}
	v(Q_i)&\geq \min_{p^k\leq j \leq q}  \left\{v\binom j {p^k} + v(P_j) + (j-p^k)v(\alpha_n)\right\} + (p^k - i)v(\alpha_n)\notag\\
	&= v(Q_{p^k}) + (p^k - i)v(\alpha_n)\notag
\end{align}
	as was to be shown.
	
	Lastly, $y_{p^k}$ is simply $v(Q_{p^k})$, which is given by (\ref{npheighteqn}).
\end{proof}

In the preceding description of the heights of the points defining $\NP_n$, one might notice that for sufficiently large $n$, the ``error terms'' $(j-p^k)v(\alpha_n)$ appearing in the minimum are very small. So we should expect the polygons $\NP_n$ to be quite similar when $n$ is large. This is the case, as we will prove shortly, although tracking these error terms make the proof less clear than we might like.

The main idea is that the height of each point defining $\NP_n$ has a main term and an error term. Sometimes, one can identify a vertex or non-vertex simply by the position of its main term relative to the other main terms, because the error is small. When vertices are not distinguished by the main term, it must be the error term distinguishing the vertex, and there is sufficient regularity in these error terms that when a vertex appears in $\NP_n$ due to the error term, it continues to do so for $\NP_{n+1}$ and so on. 

This important, but technical, geometric fact is made precise by the following lemma.
\begin{lemma}\label{verticeslemma}
	Let $m,m',m''$ and $0\leq e,e',e''\leq q-1$ be nonnegative integers, $0\leq s<t<u \leq r$ positive integers, and $|C|\leq 1$ a constant.
	
	For $n\geq 2$, define the following sequences of points:
	\begin{align*}
	\mathcal P_n &=\left(p^s, m + e \frac C {q^n}\right),\\
	\mathcal P_n' &=\left(p^t, m' + e' \frac C {q^n}\right),\\
	\mathcal P_n'' &=\left(p^u, m'' + e'' \frac C {q^n}\right).
\end{align*}
	
	Then the point $\mathcal P_n'$ lies below the line connecting the points $\mathcal P_n$ and $\mathcal P_n''$ if and only if the point $\mathcal P_{n+1}'$ lies below the line connecting the points $\mathcal P_{n+1}$ and $\mathcal P_{n+1}''$.
	
\end{lemma}
\begin{proof}
	The key point is the following: the slope of a line between any two lattice points over $p^u$ and $p^s$ has denominator $p^u - p^s$, which is always smaller than $q-1$, so if such a line doesn't pass through some lattice point, the closest it can approach that lattice point is at a vertical distance of $\frac 1 {q-1}$.
	
	With that in mind, $\mathcal P_n'$ lies below the line connecting $\mathcal P_n$ and $\mathcal P_n''$ if and only if
	\begin{equation}\label{ineqnslopeline}
	m' + e' \frac C {q^n} < \frac{p^t-p^s}{p^u-p^s} \left(m + e \frac C {q^n}\right)
	                           +\frac{p^u-p^t}{p^u-p^s} \left(m'' + e'' \frac C {q^n}\right).
\end{equation}
	Our goal is to show that (\ref{ineqnslopeline}) holds with $n+1$ in place of $n$: 
	\begin{equation}\label{goalineqnslopeline}
	m' + e' \frac C {q^{n+1}} < \frac{p^t-p^s}{p^u-p^s} \left(m + e \frac C {q^{n+1}}\right)
	                           +\frac{p^u-p^t}{p^u-p^s} \left(m'' + e'' \frac C {q^{n+1}}\right).
\end{equation}

We can see that inequality (\ref{ineqnslopeline}) roughly decomposes into two pieces: one involving only the main terms $m,m',m''$, and one involving just the error terms $e,e',e''$. This leads us to consider two cases:
	\begin{equation}\label{mainslopeineq}
	m' \leq \frac{p^t-p^s}{p^u-p^s} m
	                           +\frac{p^u-p^t}{p^u-p^s} m''
\end{equation}
	and
	\begin{equation}\label{mainslopeineqContra}
	m' > \frac{p^t-p^s}{p^u-p^s} m
	                           +\frac{p^u-p^t}{p^u-p^s} m''.
\end{equation}

\textit{Case 1.} If (\ref{mainslopeineq}) holds, then subtracting it from (\ref{ineqnslopeline}) and dividing by $q$ yields 
	\begin{equation}\label{errorslopeineq}
	e'\frac C{q^{n+1}} < \frac{p^t-p^s}{p^u-p^s} e \frac C{q^{n+1}}
	                           +\frac{p^u-p^t}{p^u-p^s} e''\frac C{q^{n+1}}.
	\end{equation}
	Adding (\ref{errorslopeineq}) back to our assumption (\ref{mainslopeineq}) yields the desired inequality (\ref{goalineqnslopeline}). These manipulations can be reversed, so (\ref{mainslopeineq}) is equivalent to (\ref{goalineqnslopeline}) in this case.
	
\textit{Case 2.} If (\ref{mainslopeineqContra}) holds instead, we will have a contradiction. By our key observation, the fact that (\ref{mainslopeineqContra}) is a strict inequality means that
	\begin{equation}\label{mainslopelowerbound}
	m' - \frac{p^t-p^s}{p^u-p^s} m - \frac{p^u-p^t}{p^u-p^s} m''\geq \frac 1 {q-1}
\end{equation}
	
	However, we can rearrange (\ref{ineqnslopeline}) to obtain
	\begin{equation}\label{ineqerrorcontra}
	m' - \frac{p^t-p^s}{p^u-p^s} m - \frac{p^u-p^t}{p^u-p^s} m'' < -e' \frac{C}{q^n} + \frac{p^t-p^s}{p^u-p^s} e \frac{C}{q^n} + \frac{p^u-p^t}{p^u-p^s} e'' \frac{C}{q^n}.
\end{equation}
	
	The left hand side is at least $\frac 1 {q-1}$ by (\ref{mainslopelowerbound}), but the right hand side is too small to allow this:
	\begin{align}
	\left|-e' \frac{C}{q^n} + \frac{p^t-p^s}{p^u-p^s} e \frac{C}{q^n} + \frac{p^u-p^t}{p^u-p^s} e'' \frac{C}{q^n}\right|
	  &= \left|-e' + \frac{p^t-p^s}{p^u-p^s} e + \frac{p^u-p^t}{p^u-p^s} e'' \right|\left| \frac{C}{q^n}\right|\notag\\
	  &\leq \left|\frac{p^t-p^s}{p^u-p^s} (q-1) + \frac{p^u-p^t}{p^u-p^s} (q-1) \right|\left| \frac{C}{q^n}\right|\notag\\
	  &= \left|q-1\right|\left|\frac{C}{q^n}\right|\notag\\
	  &\leq (q-1) \frac 1 {q^2}\notag\\
	  &< \frac 1 q.\label{finalcontra}
\end{align}

Together, (\ref{mainslopelowerbound}), (\ref{ineqerrorcontra}), and (\ref{finalcontra}) give $\frac 1 {q-1} < \frac 1 q$, clearly a contradiction.
\end{proof}

With Lemma~\ref{verticeslemma} in hand, we are ready to prove the final result of this section, a crucial input to our main results.

\begin{proposition}\label{PropVertexStability}
	There is a positive integer $V$ depending only on the polynomial $P(x)$ and the sign of $v(\alpha_0)$ such that for all $n$ sufficiently large the Newton polygon $\NP_n$ of $P(x+\alpha_n) - \alpha_{n-1}$ has exactly $V$ vertices.
	
	In fact, there are nonnegative integers $r_i, m_i, e_i$, for $1\leq i \leq V$, depending only on $P$ and $v(\alpha_0)$, and a constant $C$ which depends only on the degree $q$ and sequence of valuations $(v(\alpha_n))_{n\in\N}$, such that, for all sufficiently large $n$, the vertices of $\NP_n$ are all of the form
	$$\left(p^{r_i}, m_i + \frac{e_i}{q^n} C\right).$$
\end{proposition}
\begin{proof}
	We start by using the results of Section 2, which characterize the good behavior of ramification for large $n$, by applying Propositions~\ref{PropSemistable} and Lemma~\ref{LemmaNP}. Together, these propositions tell us that there is some $N$ such that $|v(\alpha_N)| \leq \frac 1 {q^2}$ and all the conclusions of both Proposition~\ref{PropSemistable} and Lemma~\ref{LemmaNP} hold for $n\geq N$. For the remainder of the proof, we only discuss $n\geq N$. Set $C = q^Nv(\alpha_N)$; this is independent of our choice of $N$, which we can see by again applying Proposition~\ref{PropSemistable}:
	\begin{equation}\label{Cindepeqn}
	q^nv(\alpha_n) = q^n\frac{v(\alpha_N)}{q^{n-N}} = q^N v(\alpha_N) = C,
\end{equation}
	from which it also follows that, for all $n\geq N$, $v(\alpha_n) = \dfrac C{q^n}$.
	
	Now, recall the description of $\NP_n$ given by Lemma~\ref{LemmaNP}: it is the lower convex hull of the points $(p^k,y_{p^k})$, where
	$$y_{p^k} = \min_{p^k\leq j \leq q} \left\{v\binom j {p^k} + v(P_j) + (j-p^k)v(\alpha_n)\right\}.$$
	
	Since $|v(\alpha_n)| \leq \frac 1 {q^2}$ and $|j-p^k| \leq q-1$,
	$$|(j-p^k)v(\alpha_n)| < 1,$$
	while $v\binom j {p^k} + v(P_j)$ is an integer. Moreover, all the terms $(j-p^k) v(\alpha_n)$, for $k$ fixed and $n,j$ varying, have the same sign, and so the index $j$ which achieves the minimum is determined entirely by the ``main term'' $v\binom j {p^k} + v(P_j)$ except when ties must be broken. The ties always break the same way, and depend only on the sign of $v(\alpha_0)$: in the integral case, one takes the smallest index $j$ achieving the tie value, while in the non-integral case one takes the largest such index. These are the choices which minimize the expression when there is a tie for the larger contribution of $v\binom j {p^k} + v(P_j)$.
	
	So for each $k$, the height of the point above $p^k$ is
	$$y_{p^k} = \min_{p^k\leq j \leq q} \left\{v\binom j {p^k} + v(P_j) + (j-p^k)v(\alpha_n)\right\}$$
	with the minimum achieved by a unique index $j$ between $p^k$ and $q$. Then define $M_{p^k}$ to be $v\binom j {p^k} + v(P_j)$ and $E_{p^k}$ to be $j-p^k$. The above argument shows that $j$ is independent of $n$, and hence these quantities are also independent of $n$. Clearly all are positive. Moreover, because $v(\alpha_n) = \frac{C}{q^n}$, we see that
	\begin{equation}\label{heightypkeqn}
	y_{p^k} = M_{p^k} + \frac{E_{p^k}}{q^n}C.
\end{equation}
	
	From (\ref{heightypkeqn}), what remains to be shown is that the number of vertices and the $x$-coordinates of the vertices do not depend on $n$. This follows essentially immediately from Lemma~\ref{verticeslemma} and a straightforward induction, together showing that if the Newton polygon $\NP_n$ has a vertex over $p^t$ then the Newton polygon $\NP_{n+1}$ does too, and conversely that if $\NP_n$ has no vertex over $p^t$, then neither does $\NP_{n+1}$.
	
	We know that $\NP_n$ has a vertex over $p^t$ if and only if for all $s$ and $u$ such that $s<t<u$ the point over $p^t$ lies below the line segment connecting the vertices over $p^s$ and $p^u$. If we let 
	$$m = M_{p^s},\ m' = M_{p^t},\ m'' = M_{p^u},\ e = E_{p^s},\ e' = E_{p^t},\ e'' = E_{p^u},$$
	then we are exactly in the situation to which Lemma~\ref{verticeslemma} applies: by (\ref{heightypkeqn}) the points $\mathcal P_n,\mathcal P'_n,\mathcal P''_n$ are the points over $p^s,p^t,$ and $p^u$ defining $\NP_n$, while $\mathcal P_{n+1},\mathcal P'_{n+1},\mathcal P''_{n+1}$ are the points over $p^s,p^t,$ and $p^u$ that are used to define $\NP_{n+1}$. So the lemma tells us that $\NP_n$ has a vertex over $p^t$ if and only if $\NP_{n+1}$ also has a vertex over $p^t$.
	
	Thus, by induction, all of the vertices lie over the same $x$-coordinates for all $n\geq N$, and hence their number, which we call $V$, is constant. We let $r_i$ be the exponents of the prime powers which appear as $x$-coordinates; $m_i$ be the associated main term $M_{p^{v_i}}$; $e_i$ the associated error coefficient $E_{p^{v_i}}$. The arguments above show that these do not depend on the choice of branch, only the valuations of the coefficients of $P(x)$ and the sign of $v(\alpha_0)$. We note that the subscripts indexing $m_i$ and $e_i$ are incompatible with the subscripts indexing $M_{p^k}$ and $E_{p^k}$.
		
	To conclude, we let
	$$C= \lim_{n\rightarrow \infty} q^n v(\alpha_n).$$
	As was shown in (\ref{Cindepeqn}), that sequence $q^nv(\alpha_n)$ is eventually constant, so this limit exists; clearly it only depends on $q$ and the sequence of valuations $\{v(\alpha_n)\}_{n\in\N}$. The proof above shows that $C$ plays the desired role in defining the heights of the vertices.
	
\end{proof}

\begin{definition}\label{defnlimitingramdata}
	In the notation of the preceding proposition, we define the \textit{limiting ramification data} associated to $P$ and the branch:
	\begin{align*}
	V(P,(\alpha_n)_{n\in N})&= \textrm{the number of vertices }V,\\
	R(P,(\alpha_n)_{n\in N}) &= (r_1,...,r_V),\\
	M(P,(\alpha_n)_{n\in N}) &= (m_{1},...,m_{V}),\\
	E(P,(\alpha_n)_{n\in N}) &= (e_{1},...,e_{V}),\\
	C(P,(\alpha_n)_{n\in N}) &= \textrm{the constant } C.
\end{align*}
We refer to these quantities as the ``number of vertices'', ``vertex exponents'', ``main terms'', ``error factors'', and ``error coefficient'', respectively.

\end{definition}
	Since the first vertex is over $1$ and the last vertex is $(q,0)$, defined by a minimum with just one term, we see that $r_1 = 0$ and $r_V = r$ and $m_V = e_V = 0$.
	
	As was pointed out in Proposition~\ref{PropVertexStability}, $V$, $R$, $M$, and $E$, only depend on the (ordered) valuations of the coefficients of $P$ and the sign of $v(\alpha_0)$, while $C$ depends only on the degree $q$ of $P(x)$ and the sequence $(v(\alpha_n))_{n\in N})$ of the members of the branch. The calculation of these parameters is effective, and algorithms for their computation are outlined in Section 4, along with an example.

	In fact, the only ineffective step in our results occurs in Proposition~\ref{PropSemistable} -- the proof of (b) and (c) does not give an effective determination of ``sufficiently large''. There are some cases where this can be circumvented; for instance, if $v(\alpha_0) = 1$ then it is straightforward to see that, for all $n$, $P(x)-\alpha_n$ is Eisenstein, which implies (b) and (c) hold for all $n$. More generally, it follows from our proof of Proposition~\ref{PropSemistable} that if there is some $N$ such that $v(\alpha_N)$ is not divisible by $p$ and has smaller valuation than any coefficient of $P(x)$, then (b) and (c) hold for all $n\geq N$.

One can see quite readily from Proposition~\ref{PropVertexStability} that the polygons $\NP_n$ have a pointwise limit (viewing them as functions on $\mathbb R_\geq 0$). Some of what follows can be described in terms of that limiting polygon, and at times more simply -- for instance, one could avoid using Lemma~\ref{verticeslemma}. However, valuable information is lost when working with this limit polygon: it may have fewer vertices than the actual Newton polygons $\NP_n$ (this occurs when main terms of vertices, $(p^{r_i},m_i)$, are collinear). The number of vertices $V$ is extremely important for our main result and applications, because $V-1$ is the slope of the linear change of index in our main result. Additionally, it is appealing to have such an exact description of $\NP_n$.

\section{Main Results}

As mentioned in the introduction, our second main result describes the higher ramification after extending the ground field and adjusting the index. This is because the dynamics of the ramification can take some time to stabilize. As we will see shortly, the previous section amounted to showing that the ramification actually does stabilize. So for technical simplicity, we will now introduce some assumptions to the effect that we have already reached the region of stable behavior (in other words, that the results of the preceeding section hold immediately for $P$ and $\alpha_0$, without first replacing the base point by some $\alpha_N$). At the end we will explicitly work out the reduction of the general case to the stable case, and the adjustments required. This leads us to introduce the following property:
\begin{definition}
	A pair $(P,\alpha_0)$ satisfies (H) if they satisfy the \textit{conclusions} of Propositions~\ref{PropSemistable} and~\ref{PropVertexStability} for all $n$, without the qualification ``for sufficiently large $n$''. 
\end{definition}

And so Propositions~\ref{PropSemistable} and~\ref{PropVertexStability} tell us that even if $P$ and $\alpha_0$ do not satisfy (H), there is some $N$ such that $P$ and $\alpha_N$ do. In the Galois case, this is equivalent to replacing an (infinite) profinite group with a finite-index subgroup which, hopefully, retains a lot of information about the original group.

Besides this, it is also necessary to introduce a ``tameness'' assumption, that $d$ is not divisible by $p$. Recall that $d=\lim_{n\rightarrow\infty} d_n$ is the eventual valuation of $\alpha_n$ with respect to a valuation that sends $\pi_n$ to $1$, and that this limit exists was shown in Proposition~\ref{PropSemistable}(c). In what follows, we will want to take a $d$th root of the unit $u_n$ inside $K_n$. Recall that the unit $u_n$ was defined by $\alpha_n = u_n\pi_n^{d_n}$, and so the presence of this $d$th root allows us to take a different choice of uniformizer $\pi_n$, such that $\alpha_n = \pi_n^d$. This $d$th root is not necessarily in $K_n$, but if $p$ does not divide $d$, then we can obtain a $d$th root of $u_n$ after an unramified extension of $K_n$, which does not change the ramification along the branch. However, if $p$ divides $d$ then the $d$th root of $u_n$ may only appear in a ramified extension of $K_n$, and this extra ramification interferes with our ability to extract information about ramification prior to including the $d$th root. We hope that this restriction can be relaxed in some or all cases -- the study of some special cases suggests that if $d=d_0 p^m$ where $p\nmid d_0$, then our results still hold with $d_0$ in place of $d$. An unfortunate downside of this restriction is that it means our results are not base-change invariant -- if we replace $K$ by an extension with ramification index divisible by $p$ and linearly disjoint from $K_\infty$, then $p$ is guaranteed divide $d$. Luckily, we at least have invariance under \textit{tame} base change.

Given these assumptions, our next goal is verify that the extension $K_\infty/K$ is arithmetically profinite and compute its Hasse-Herbrand function, under (H) and the tameness assumption. We will break up the computation of the Hasse-Herbrand function of $K_\infty/K$ into calculating the Hasse-Herbrand functions for the intermediate extensions $K_n/K_{n-1}$, composing those functions to obtain the Hasse-Herbrand function of $K_n/K$, and then pass to the limit. As mentioned in the introduction, we avoid assuming any of our extensions are Galois (indeed, one would expect this to be rare in general) so to study higher ramification, we employ the techniques explained by Lubin~\cite{Lubin}. The reader is advised to take some care in passing between this and other sources (such as Serre~\cite{Serre}) since the ramification groups may be numbered differently; we adopt Lubin's convention.

For convenience, we remind the reader of two important polygons: the \textit{level $n$ Newton polygon} $\NP_n$ and its dual, the \textit{level $n$ Newton copolygon} denoted $\co\NP_n$. The former is the Newton polygon of $P(x+\alpha_n) - \alpha_{n-1}$, while the latter is its dual, meaning that $\co\NP_n$ has a vertex whose $x$-coordinate is the negative of that slope, and the slopes of $\co\NP_n$ slopes are the values $p^{r_i}$, in decreasing order. As such, the copolygon $\co\NP_n$ has one fewer vertex  than the polygon $\NP_n$. The assumption (H) amounts to the following explicit description of $\NP_n$: it is the lower convex hull of the following points determined by the limiting ramification data:
	$$\left(p^{r_i}, m_i + e_i \frac{C}{q^n}\right)\ \ \ \ 1\leq i \leq V.$$

\begin{proposition}\label{smallHH}
	Suppose the pair $(P,\alpha_0)$ satisfies (H) and that $p$ does not divide $d$. Then the graph of the Hasse-Herbrand transition function $\phi_n$ for $K_n/K_{n-1}$ relative to $K$ can be obtained by applying the following three transformations to the copolygon $\co\NP_n$:
	\begin{enumerate}[(1)]
	\item Increase the $x$-coordinates of each vertex by $\mathrm{sgn}(v(\alpha_0))(d-1)v(\pi_n)$, while modifying the $y$-coordinates to preserve the slopes of the segments between them.
	\item Stretch horizontally by a factor of $e_{K/E}q^n$.
	\item Stretch vertically by a factor of $e_{K/E}q^{n-1}$.
\end{enumerate}
	
	The first slope of $\phi_n$ is $1$ and the last slope of $\phi_n$ is $1/q$. The $x$-coordinates of the first and last vertices of $\phi_n$, are respectively,
	$$- e_{K/E}q^n(\mathrm{shallowest\ slope\ of\ } \NP_n) + \mathrm{sgn}(v(\alpha_0))(d-1)v(\alpha_0)$$
	and
	$$- e_{K/E}q^n(\mathrm{steepest\ slope\ of\ } \NP_n) + \mathrm{sgn}(v(\alpha_0))(d-1)v(\alpha_0).$$
\end{proposition}
\begin{proof}
	We will prove this in full for the integral case, where $v(\alpha_0) > 0$ and hence $d \geq 1$, and at the end indicate the minor modifications necessary for the non-integral case.
	
	Let $f(x)$ be the minimal polynomial for $\pi_n$ over $K_{n-1}$. The Hasse-Herbrand function for $K_n/K_{n-1}$ can be obtained by applying stretches (2) and (3) to the Newton copolygon of $f(x+\pi_n)$~\cite[Definition 5]{Lubin}. In Lubin's notation, we are taking $K=K_n$, $k=K_{n-1}$, and $k_0 = E$, and $\Psi_{v,F}$ is the copolygon of $f(x+\pi_n)$; the claimed scaling factors are obtained by expanding $e_{K/k_0} = e_{K_n/K}e_{K/E} = q^n e_{K/E}$ and similarly $e_{k/k_0} = q^{n-1}e_{K/E}$. So we only need to show that the copolygon of $f(x+\pi_n)$ can itself be obtained by applying (1) to $\co\NP_n$.
	 
	 In terms of Newton polygons, (1) is equivalent to \textit{decreasing} all of the slopes of $\NP_n$ by $(d-1)v(\pi_n)$ (there is a sign change in the duality between polygon and copolygon!). The Newton polygons of $f(x+\pi_n)$ and $P(x+\alpha_n) - \alpha_{n-1}$ encode the valuations of the roots of the corresponding polynomials. These roots are of the form $\pi_n^\sigma - \pi_n$ and $\alpha_n^\sigma - \alpha_n$, respectively, for $\sigma \in \Gamma_K$, and our task is to relate their valuations.
	 
	 In the integral case, we want to show that, for all $\sigma \in \Gamma_K$,
	 $$v(\pi_n^\sigma - \pi_n) = v(\alpha_n^\sigma - \alpha_n) - (d-1)v(\pi_n).$$
	 
	 Recall that we selected uniformizers $\pi_n$ and units $u_n$ such that $\alpha_n = u_n \pi_n^{d_n}$. By (H), $d_n = d$ does not vary with $n$, and we also assumed it is not divisible by $p$. As such, $u_n$ admits a $d$th root after at most an unramified extension; the transition function is insensitive to base change by unramified extensions. In other words, extending $K_{n-1}$ and $K_n$ by their unramified extension of degree $d$ guarantees the presence of $\sqrt[d]{u_n}$ in our field without affecting the transition function. So altering our choice of $\pi_n$, we may write $\alpha_n = \pi_n^d$. This allows us to compare the valuations more directly:
	\begin{align}
	\alpha_n^\sigma - \alpha_n &= (\pi_n^\sigma)^d - \pi_n^d\label{prodeqnpiline}\\
	    &= \prod_{\zeta^d = 1} (\pi_n^\sigma - \zeta \pi_n)\label{prodeqnpi}
\end{align}
	
	Of the terms in the product (\ref{prodeqnpi}), we are only interested in $v(\pi_n^\sigma - \pi_n)$. To manage the others, notice that
	\begin{equation}\label{valunityeqn}
		v(\pi_n^\sigma - \zeta \pi_n) = v(\pi_n) + v\left(\frac{\pi_n^\sigma}{\pi_n} - \zeta\right).
\end{equation}
	
	If $v\left(\frac{\pi_n^\sigma}{\pi_n} - \zeta\right)$ is positive, then $\frac{\pi_n^\sigma}{\pi_n}$ is necessarily a $d$th root of unity modulo $\pi_n$. On the other hand, the norm from $K_n$ to $K$ of $\frac{\pi_n^\sigma}{\pi_n}$ is just $1$; but viewed in the residue field, the norm is just the $q$th power. Therefore, in the residue field, $\frac{\pi_n^\sigma}{\pi_n}$ is both a $d$th root of unity and a $q$th root of unity. Because $p\nmid d$, this is only possible if $\zeta = 1$. In all other cases, $v\left(\frac{\pi_n^\sigma}{\pi_n} - \zeta\right) = 0$. Thus, (\ref{valunityeqn}) simplifies to just $v(\pi_n)$ whenever $\zeta\neq 1$, and so the valuation of (\ref{prodeqnpi}) becomes
	$$v(\alpha_n^\sigma - \alpha_n) = v(\pi_n^\sigma - \pi_n) + (d-1)v(\pi_n)$$
	or equivalently 
	$$v(\pi_n^\sigma - \pi_n) = v(\alpha_n^\sigma - \alpha_n) - (d-1)v(\pi_n),$$
	which is exactly the statement to which we reduced the main part of this proposition for the integral case.
	
	For the non-integral case, when $d$ is negative, we must instead work with
	$$\frac 1 {\alpha_n^\sigma} - \frac 1 {\alpha_n} = (\pi_n^\sigma)^{|d|} - \pi_n^{|d|}.$$
	The left hand side can be written as
	$$\frac{\alpha_n-\alpha_n^\sigma}{\alpha_n\alpha_n^\sigma}$$
	which has valuation
	$$v(\alpha_n-\alpha_n^\sigma) - 2v(\alpha_n).$$
	
	Recall too that $v(\alpha_n) = dv(\pi_n)$. After replacing the left hand side of (\ref{prodeqnpiline}) and rearranging to move the $2dv(\pi_n)$ to the right hand side, the remainder of the argument proceeds essentially unchanged until the end, where incorporating the extra $2dv(\pi_n)$ gives rise to the $\sgn(v(\alpha_0))$ in the statement of the proposition.
			
	Finally, by inspecting the transformation of $\co\NP_n$ into $\phi_n$, one can see that the first and last slopes of $\phi_n$ are $\frac {e_{K/E}q^{n-1}} {e_{K/E}q^n} = \frac 1 q$ multiplied by the first and last slopes of $\co\NP_n$. The first and last slopes of $\co\NP_n$ are the first and last $x$-coordinates of vertices of $\NP_n$, which are $1$ and $q$, so together we see that the first and last slopes of $\phi_n$ are $1$ and $\frac 1 q$, as claimed. Likewise, the $x$-coordinates can be obtained from the duality of $\co\NP_n$, which turns negative slopes of $\NP_n$ into $x$-coordinates of vertices, then modified according to the first two transformations.
\end{proof}
\begin{remark}
	We use the assumption $p\nmid d$ in two places: to take a $d$th root of $u_n$, and that the $d$th roots of unity are distinct modulo $p$ to control $v(\pi_n - \zeta \pi_n)$. The former seems to be the greater obstacle to characterizing ramification in the general case.
\end{remark}

The essence of the preceding proposition is that the ramification-theoretic properties of these extensions are somewhat stable. Neglecting scaling, all the Hasse-Herbrand functions look like a small shift of $\co\NP_n$, and the copolygon itself changes little as a function of $n$, in a way which is described very precisely by Proposition~\ref{PropVertexStability}.

The most difficult step would appear to be composing the intermediate Hasse-Herbrand functions $\phi_1,\phi_2,...,\phi_n$ to obtain the Hasse-Herbrand function $\Phi_n$ for $K_n/K$. However, this is straightforward if we can ensure that the $\phi_n$'s behave sufficiently well. Since $\phi_n$ is the identity along its first segment, one might hope that the domain on which it is the identity includes all of the vertices of $\Phi_{n-1}$.

Unfortunately, this is too much to expect in general, but we can give a characterization of when these functions do have sufficiently large identity segments in terms of the limiting copolygon. We can show that both post-critically bounded polynomials (of the appropriate form) and polynomials of prime degree exhibit this good behavior with respect to composition of the above Hasse-Herbrand functions.

\begin{proposition}\label{PropComposition}
	Suppose $(P,\alpha_0)$ satisfy (H) and that $p$ does not divide $d$.
	
	For $n\geq 2$, the leftmost vertex of $\phi_n$ has strictly larger $x$-coordinate than that of the rightmost vertex of $\phi_{n-1}$ if the limiting Newton polygon has just one slope, or if
	 $$-q\frac{m_{V} - m_{V-1}}{p^{r_V} - p^{r_{V-1}}} > -\frac{m_2 - m_1}{p^{r_2} - p^{r_1}} + \frac{2}{p-1}|v(\alpha_{0})|.$$
	 
	 The various $m_i$ and $r_i$ are the quantities given by the limiting ramification data of Definition \ref{defnlimitingramdata}.
\end{proposition}
\begin{proof}
	By the final statement of Proposition~\ref{smallHH}, we can rewrite the claim about the $x$-coordinates of those vertices in terms of the slopes of $\NP_n$ and $\NP_{n-1}$. We want
	$$- e_{K/E}q^n(\mathrm{shallowest\ slope\ of\ } \NP_n) + \sgn(v(\alpha_0))(d-1)v(\alpha_0)$$
	to be strictly larger than
	$$- e_{K/E}q^{n-1}(\mathrm{steepest\ slope\ of\ } \NP_{n-1}) + \sgn(v(\alpha_0))(d-1)v(\alpha_0).$$
	
	For convenience, let's name the negatives of these slopes: let 
	$$s = -(\mathrm{shallowest\ slope\ of\ } \NP_n)\ \ \ \ s'= -(\mathrm{steepest\ slope\ of\ } \NP_{n-1}).$$
	
	Now we can simplify and rewrite the target inequality as
	$$q s > s'.$$
	 When there is just one slope, $s=s'$ and the inequality obviously holds. Otherwise, there are two slopes.
	 
	Now, the height of the vertex over $p^{r_i}$ is given by $m_i + e_i \frac {v(\alpha_0)}{q^n}$. So in expressing the slopes of the segments between our vertices of interest in these terms, the quantities
	$$t= -\frac{m_{V} - m_{V-1}}{p^{r_V} - p^{r_{V-1}}}\ \ \ \ \textrm{and}\ \ \ \ t'=-\frac{m_2 - m_1}{p^{r_2} - p^{r_1}}$$
	in the statement of the proposition are the (negative) contributions of the ``main terms''  to the slopes $s$ and $s'$. In light of this interpretation, we can write
	\begin{align*}
	s-t   &= \frac{e_{V-1} - e_{V}}{p^{r_B} - p^{r_{V-1}}}\frac{v(\alpha_0)}{q^n},\\
	s'-t' &= \frac{e_1 - e_2}{p^{r_2} - p^{r_1}}\frac{v(\alpha_0)}{q^{n-1}}.
	\end{align*}
	As was remarked previously, $r_1 = 0$, $r_V = r$ and $e_V = 0$, because the first vertex lies over $1$, while the last vertex is $(q,0)$. 
	
	To summarize, the hypothesis of the proposition is
	$$qt > t' + \frac 2 {p-1} v(\alpha_0),$$
	and we have some $s,s'$ such that
	\begin{align*}
	s-t   &= \frac{e_{V-1}}{q - p^{r_{V-1}}}\frac{v(\alpha_0)}{q^n},\\
	s'-t' &= \frac{e_1 - e_2}{p^{r_2} - 1}\frac{v(\alpha_0)}{q^{n-1}},
\end{align*}
	and our goal is
	$$qs > s'.$$
	
	But then it is enough for our two errors $q(s-t)$ and $s'-t'$ to be small enough that their sum is less than $\frac 2 {p-1} |v(\alpha_0)|$ in absolute value, as then adding these error terms to the inequality we initially assumed will preserve the inequality up to the loss of margin of error, $\frac 2 {p-1} |v(\alpha_0)|$, that we allowed ourselves. To prove that the sum of $q(s-t)$ and $s'-t'$ is small enough, it suffices to show that each is at most $\frac {|v(\alpha_0)|} {p-1}$. And indeed:
	\begin{align*}
	|q(s-t)| &= q\frac{e_{V-1}}{q - p^{r_{V-1}}}\frac{|v(\alpha_0)|}{q^n}\\
	         &\leq q\frac{q-1}{q - p^{r-1}}\frac{|v(\alpha_0)|}{q^n}\\
	         &< \frac{1}{q-p^{r-1}}\frac{|v(\alpha_0)|}{q^{n-2}}\\
	         &\leq \frac{|v(\alpha_0)|}{p-1},
	\end{align*}
	and
	\begin{align*}
	|s'-t'| &= \frac{|e_1 - e_2|}{p^{v_2} - 1}\frac{|v(\alpha_0)|}{q^{n-1}}\\
	         &\leq \frac{q-1}{p^{v_2} - 1}\frac{|v(\alpha_0)|}{q^{n-1}}\\
	         &\leq \frac{q-1}{p - 1}\frac{|v(\alpha_0)|}{q^{n-1}}\\
	         &< \frac{1}{p - 1}\frac{|v(\alpha_0)|}{q^{n-2}}\\
	         &< \frac{|v(\alpha_0)|}{p - 1},
	\end{align*}
	where, on the second line, we use $|e_1-e_2| \leq q-1$ rather than $\leq 2(q-1)$ because we know that $e_1$ and $e_2$ are both nonnegative. Both inequalities also require that $n\geq 2$ so that $\frac 1 {q^{n-2}}$ is at most $1$. 
\end{proof}
\begin{corollary}\label{CorCompp}
	Assume (H) and that $p\nmid d$. If $P$ has degree $q=p$, then it satisfies Proposition~\ref{PropComposition}.
\end{corollary}
\begin{proof}
	Immediate, as in this case the limiting Newton polygon can only have vertices over $1$ and $p$, hence it has just a single slope.
\end{proof}

\begin{corollary}\label{CorCompPCB}
	Assume (H) and that $p\nmid d$. If $P$ is post-critically bounded and
	$$|v(\alpha_0)| < \frac{p-1}{2}v(p),$$
	then the pair $(P,\alpha_0)$ satisfies the hypotheses of Proposition~\ref{PropComposition}. 
\end{corollary}
\begin{proof}
	We will verify directly that Proposition~\ref{PropComposition} applies. If the limiting Newton polygon has just one slope, we are done. Otherwise, assume it has at least two. Then we want to verify that the following inequality holds:
	 \begin{equation}-q\frac{m_{V} - m_{V-1}}{p^{r_V} - p^{r_{V-1}}} > -\frac{m_2 - m_1}{p^{r_2} - p^{r_1}} + \frac{2}{p-1}|v(\alpha_{0})|. \label{corPCBineq}\end{equation}
	 
	 Recall Proposition~\ref{PropPCB}, which says that $\frac{P'(x)}{q}$ has integral coefficients. The first vertex of $\NP_n$ is $(1,v(P'(\alpha_n))$, and so its height is at least $v(q) = rv(p)$. Moreover, from Lemma~\ref{LemmaBinomialDrop}(ii), we know that the height drop between vertices over $p^s$ and $p^u$ is at most $(u-s)v(p)$; in our notation, $-(m_i - m_j) \leq (r_i-r_j)v(p)$ for $i\geq j$. Recall as well that $r_1 = 0$ and $r_V = r$.

	 Working with the right hand side of (\ref{corPCBineq}), this means that
	 \begin{align}
	 	-\frac{m_2 - m_1}{p^{r_2} - p^{r_1}} + \frac 2{p-1}|v(\alpha_0)|\notag
	 	    &=-\frac{m_2 - m_1}{p^{r_2} - 1} + \frac 2{p-1}|v(\alpha_0)|\notag\\
	 	    &\leq\frac{r_2 v(p)}{p^{r_2} - 1} + \frac 2{p-1}|v(\alpha_0)|\notag\\
	 	    &<\frac{v(p)}{p - 1} + \frac 2{p-1}|v(\alpha_0)|\notag\\
	 	    &<\frac{v(p)}{p - 1} + v(p)\notag\\
	 	    &=\frac p {p-1}v(p).\label{corPCBright}
	 \end{align}
	 We use our assumption about $|v(\alpha_0)|$ on the second to last line.

	Now, let us treat the left hand side of (\ref{corPCBineq}). Recall $r_V = r$, $1\leq r_{V-1} \leq r-1$, $m_1 \geq rv(p)$ and $m_V = 0$. From this it follows that 
	$$-(m_{V} - m_{V-1}) = m_{V-1} \geq (r - r_{V-1})v(p)$$
	
	Applying this to the left hand side of (\ref{corPCBineq}), we see:
	
	\begin{align}
		-q\frac{m_{V} - m_{V-1}}{p^{r_V} - p^{r_{V-1}}} \notag
			&\geq  p^r \frac {(r - r_{V-1})v(p)} {p^r - p^{r_{V-1}}} \\
			&\geq  \frac {r - r_{V-1}} {1 - p^{r_{V-1} - r}} v(p)\label{corPCBtricky}\\
			&\geq  \frac p {p-1}v(p).\label{corPCBleft}
	\end{align}
	
	Going from (\ref{corPCBtricky}) to (\ref{corPCBleft}) is slightly tricky; if $r_{V-1} = r-1$ then the two are equal, while if $r_{V-1} < 1$ then in the fraction term of (\ref{corPCBtricky}), the numerator $(r - r_{V-1})$ is at least $2$ and the denominator $(1 - p^{r_{V-1} - r})$ is at most $1$, hence the whole quantity is at least $2v(p) \geq \frac p{p-1}v(p)$ (sharp when $p=2$).

	Clearly, (\ref{corPCBright}) and (\ref{corPCBleft}) yield the desired inequality (\ref{corPCBineq}).

\end{proof}

\begin{remark}
	Notably, when $p$ is at least $5$, the inequality in the proposition is always satisfied for $v(\alpha_0) = 1$.
\end{remark}

It still remains to compose our Hasse-Herbrand functions. The conclusion of Proposition~\ref{PropComposition} describes the ``good behavior'' that we want in order for the Hasse-Herbrand functions to compose well: the first vertex of $\phi_n$ should have larger $x$-coordinate than the last vertex of $\phi_{n-1}$. When this happens, the higher ramification behavior of the branch is quite well-controlled, and highly regular. From working with explicit examples, it is clear that this happens in many situations besides those described by Proposition~\ref{PropComposition} or Corollaries~\ref{CorCompp} and~\ref{CorCompPCB}. This leads us to introduce the following definition:
\begin{definition}
	A branch associated to $P$ and $\alpha_0$ over $K$ is said to be \textit{tamely ramification-stable} if $p\nmid d$, and the pair satisfies (H) and the conclusions of Propositions~\ref{smallHH} and~\ref{PropComposition}.
	
	A branch is said to be \textit{potentially} tamely ramification-stable if there is some $N$ such that upon replacing $K$ by $K_N$ and re-indexing the branch to be based at $\alpha_N$ it is tamely ramification-stable.
\end{definition}
\begin{remark}
	In our definition, ``tamely'' refers to the restriction that $p\nmid d$. We expect that even if $p|d$, such branch extensions would still exhibit this kind of ramification stability. However, the precise expressions given in Proposition~\ref{smallHH}, particularly the $(d-1)v(\pi_n)$ term, may not correctly describe these cases.
\end{remark}
\begin{proposition}\label{PropPotentiallyTRS}
	Suppose that $p\nmid d$. If $P(x)$ has prime degree or is post-critically bounded, then any branch associated to $P(x)$ is potentially tamely ramification-stable.
\end{proposition}
\begin{proof}
	Propositions~\ref{PropSemistable} and~\ref{PropVertexStability} ensure that for all sufficiently large $N$, (H) is satisfied whenever $K$ is replaced by $K_N$ and the branch is modified to start at $\alpha_N$. 
	
	For polynomials of prime degree and post-critically bounded polynomials, Corollaries~\ref{CorCompp} and~\ref{CorCompPCB}, respectively, prove that any branch satisfies the conclusion of Proposition~\ref{PropComposition} after possibly increasing $N$.
\end{proof}

From the proof of Proposition~\ref{PropComposition}, we know that if $p\nmid d$, a branch is potentially tamely ramification stable when, roughly, the first (steepest) slope of $\NP_{n-1}$ is not more than $q$ times steeper than the last (shallowest) slope of $\NP_n$. This property depends only on $P(x)$, not the branch. For it to fail, the first vertex of $\NP_n$ must be relatively high compared to the others, which seems unlikely based on the structure of the minima that describe the heights of these vertices.

Before proceeding, recall the following definition:
\begin{definition}[\cite{Lubin}]
	The \textit{altitude} of an extension $E/K$ with transition function $\Psi(x)$ is the height of the rightmost vertex of $\Psi(x)$; at times we may abbreviate this as the altitude of $\Psi(x)$.
\end{definition}

\begin{proposition}\label{mediumHH}
	Suppose our branch, associated to $(P,\alpha_0)$, is tamely ramification-stable over $K$. Let $V$ be the number of vertices from the limiting ramification data.
	
	Then the Hasse-Herbrand transition function $\Phi_n(x)$ for $K_n/K$ is a piecewise linear function which satisfies the following properties:
	\begin{enumerate}[\indent 1.]
	\item $\Phi_n(x)$ has $(V-1)n$ vertices,
	\item the last (rightmost) vertex of $\Phi_n(x)$ has the same $x$-coordinate of the last vertex of $\phi_n$,
	\item the final slope (of the ray extending rightward from the last vertex) of $\Phi_n(x)$ is $1/q^n$,
	\item $\Phi_n(x)$ coincides with $\Phi_{n-1}(x)$ for $x$ smaller than the last coordinate of $\Phi_{n-1}$,
	\item the altitude of $\Phi_n(x)$ is strictly greater than the altitude of $\Phi_{n-1}$ and is unbounded as a function of $n$.
\end{enumerate}
\end{proposition}
\begin{proof}
	By transitivity, $\Phi_n(x) = \Phi_{n-1} \circ \phi_n (x)$, so it is natural to proceed by induction. The base case is $\Phi_1 = \phi_1$, where there is nothing to prove: the shape of this function has been described explicitly already and satisfies all of the above conditions.
	
	The first vertex of $\phi_n(x)$ is after the last vertex of $\Phi_{n-1}(x)$, and $\phi_n(x)$ is the identity up to its first vertex, so property (4) follows. The $x$-coordinate of the last vertex of $\phi_n(x)$ is after that of the last vertex of $\Phi_{n-1}(x)$, and so it remains the $x$-coordinate of the last vertex of $\Phi_n(x)$, verifying property (2). Moreover, after that point, we add $V-1$ new vertices, from those of $\phi_n$, yielding (1). The final segment of $\phi_n$ corresponds to the final vertex $(q,0)$ of $\NP_n$, and hence has slope $\frac 1 q$. By inspection, the final slope of $\Phi_n$ is the product of the final slope of $\Phi_{n-1}$, which is $\frac{1}{q^{n-1}}$ and the final slope of $\phi_n$, which is $\frac 1 q$, so together the final slope is $\frac 1 {q^n}$, which is (3).
	
	Finally, the altitude is the height of the last vertex of $\Phi_n(x)$, which lies over the last vertex of $\phi_n$. By Proposition~\ref{smallHH} combined with the limiting ramification data, we can express the $x$-coordinates of the last vertices of $\Phi_{n}$ and $\Phi_{n-1}$ as
	$$Aq^n + B\ \ \textrm{and}\ \  Aq^{n-1} + B,$$
	respectively, where $A$ and $B$ are positive constants which do not depend on $n$. The constant $A$ comes from the part of the slope associated to the main terms, while $B$ comes from the error terms plus the shift by $\mathrm{sgn}(v(\alpha_0))(d-1)v(\pi_n)$, and both incorporate the scaling by $e_{K/E}$.
	
	Between these two vertices, the slopes of $\Phi_n(x)$ must be at least $\frac p {q^n}$ because the last (shallowest) slope is $\frac 1 {q^n}$ and the slopes are all powers of $p$. Then we can estimate the difference in altitudes as follows
	\begin{align*}
	\textrm{altitude}(\Phi_n) - \textrm{altitude}(\Phi_{n-1})
	    &\geq \frac p {q^n} (Aq^n + B - (Aq^{n-1} + B))\\
	    &\geq A \left(p - \frac p q\right).
\end{align*}
	Thus the gap between consecutive altitudes is bounded below by a positive constant which does not depend on $n$, and so the altitudes are unbounded as $n$ grows.
\end{proof}

With this setup, our main theorem falls readily into place:

\TheoremMain
\begin{proof}
	We first show that $K_\infty/K$ is arithmetically profinite. As explained in Wintenberger~\cite{Wintenberger}, we simply need a filtration of elementary extensions whose altitudes tend to infinity. Because $\Phi_n$ restricts to $\Phi_{n-1}$, the elementary subextensions of $K_n$ inside $K_{n-1}$ are all of the elementary subextensions of $K_{n-1}$, which gives us our tower. The altitude of $K_n$ tends to infinity by Proposition~\ref{mediumHH}, hence the heights of these elementary subextensions do as well. From this we see that the extension is arithmetically profinite, and that its Hasse-Herbrand function $\Phi(x)$ is given by the pointwise limit of the intermediate Hasse-Herbrand functions $\Phi_n(x)$. Further, by Proposition~\ref{mediumHH}, $\Phi_{n}(x)$ coincides with $\Phi_{n-1}(x)$ up to the last vertex of $\Phi_{n-1}(x)$, and so the same holds for $\Phi(x)$: whenever $x$ is smaller than the $x$-coordinate of the last vertex of $\Phi_{n}(x)$, we have $\Phi(x) = \Phi_{n}(x)$.
	
	The altitude of $K_n$ over $K$ is the same as the height of the $(V-1)n$th vertex of $\Phi$, again by our assumption that the branch is tamely ramification-stable. That altitude is strictly less than the height of the $((V-1)n+1)$th vertex of $\Phi$, and so $K_n \subseteq K_\infty^{((V-1)n + 1)}$. On the other hand, the final slope of $\Phi_n(x)$ is $\frac 1 {q^n}$, by Proposition \ref{mediumHH}. Since $\Phi(x) = \Phi_n(x)$ up to the the $((V-1)N+1)$th vertex, this is the same as the slope of $\Phi(x)$ going into the $((V-1)N+1)$th vertex, so the degree of $K_\infty^{((V-1)n + 1)}$ over $K$ is $q^n$, which is the same as the degree of $K_n$ over $K$. Thus the two fields are equal, as claimed.
\end{proof}

\CorollaryMain
\begin{proof}
	If $P(x)$ has prime degree or is post-critically bounded and $p\nmid d$, then any nontrivial branch associated to it is potentially tamely ramification-stable by Proposition~\ref{PropPotentiallyTRS}. Recall that this means that there is an $N$ such that after restricting our branch to start at $\alpha_N$ it is tamely ramification-stable over $K_N$.
	
	To keep our indexing clear, set $\beta_n = \alpha_{N+n}$, $L= K(\alpha_N)$, $L_n = L(\beta_n)$, and $L_\infty = \bigcup L_n$. Clearly $L_\infty = K_\infty$. Then our main result, Theorem~\ref{TheoremMain}, applies to this branch, and so $K_\infty=L_\infty$ as an extension of $L=K_N$ is arithmetically profinite and
	$$L_n = L_\infty^{((V-1)n + 1)}.$$
	
	Translating from $L$ to $K$, we see that $K_n = L_{n-N}$ if $n\geq N$. So making this change of index, we see that
	$$K_n = L_{n-N} = L_\infty^{((V-1)(n-N) + 1)}$$
	for $n\geq N$, for the upper numbering relative to $L=K_N$. When $n < N$, $((V-1)(n-N)+1)$ is negative, which is handled by our convention for negative-indexed elementary subfields, that they are simply the ground field. Thus replacing $K$ by $K_N$ yields the claimed statement for all $n$.
\end{proof}

\section{Applications and Effectivity}
\subsection{A question of Berger.}
As our first application, we can offer a partial answer to a question raised by Berger~\cite{BergerIterated}. That paper considers extensions of the same type studied here, though with two restrictions: the degree is the size of the residue field, and the base point is a uniformizer. An important intermediate result of that paper is the implication
$$K_\infty/K\ \textrm{ Galois}\ \ \Rightarrow\ \ K_\infty/K\ \textrm{ abelian}.$$
Berger asks if there is a more direct or elementary proof of this fact: the two proofs we are aware of, due to Berger~\cite{BergerLifting} and Cais-Davis~\cite{CD}, use quite sophisticated machinery. Our results allow us to give such an elementary proof in some cases.

Let us outline Berger's use of this fact: if $K_\infty/K$ is abelian, then $K_n/K$ is also abelian, and in particular normal. When $K_n/K$ is normal and the degree of $K_n/K_{n-1}$ is $q$, one can define, for each $\sigma \in \Gamma_K$, a power series $\Col_\sigma \in K[[T]]$ such that $\Col_\sigma(0) = 0$ and $\Col_\sigma(\alpha_n) = \alpha_n^\sigma$ (generalized Coleman power series). This power series commutes with $P$, and so by a result of Lubin~\cite{LubinDynamical}, that power series is determined by the coefficient of its linear term, which gives a character from $\Gamma_K$ to $\OK^*$. This character is injective, because the action on the branch determines the action everywhere in the extension, since the branch generates the extension. Berger then goes on to study this character in detail.

But the logic flows the other way too: if we know that $K_n/K$ is normal for some other reason, then we can construct these power series and the associated injective character, which would prove that $K_\infty/K$ is abelian. And indeed, the elementary subfields of $K_\infty$ over $K$ are all normal over $K$ if $K_\infty/K$ is normal. Thus if one were to know that for all $n$ there exists an $m$ such that $K_n = K_\infty^{(m)}$ for some $m$, as in our main theorem, then $K_\infty/K$ must be abelian.

\ThmBergerApp
\begin{proof}
	Because $\alpha_0$ is a uniformizer, all of the polynomials $P^n(x) - \alpha_0$ are Eisenstein, so they are irreducible and give rise to a totally ramified extension of degree $q^n$. This means that $d=1$ and that $[K_n:K_{n-1}] = q$ for all $n$.
	
	The branch is tamely ramification stable, so we may apply Theorem~\ref{TheoremMain}, to conclude that for all $n$, the extension $K_n/K$ is elementary, and therefore also Galois because $K_\infty/K$ is Galois.
	
	Now let $\sigma \in \Gal(K_\infty/K)$. Because $K_n/K$ is normal, $\alpha_n^\sigma$ is in $K_n = K(\alpha_n)$. The sequence $(\alpha_n^\sigma)_{n\in \N}$ is itself a branch, and by our assumption that $p$ is odd and the irreducibility of $P(x) - \alpha_{n-1}$, we see that $N^{K_n}_{K_{n-1}} (\alpha_n) = \alpha_{n-1}$. This means that we can use Berger's construction (Theorem~3.1~\cite{BergerIterated}) to produce a uniquely determined series $\Col_\sigma \in \OK[[T]]$ which acts by $\Col_\sigma(\alpha_n) = \alpha_n^\sigma$ and commutes with $P(x)$. This gives rise to a character $\chi$ from $\Gal(K_\infty/K)$ to $\OK^*$ given by $\chi(\sigma) = \Col_\sigma'(0).$
	
	Since $\Col_\sigma$ commutes with $P(x)$ and $P'(0)$ is neither zero nor a root of unity, the series $\Col_\sigma$ is determined by $\Col_\sigma'(0)$ by Proposition 1.1 of Lubin~\cite{LubinDynamical}. Since $\Col_\sigma$ also determines the action of $\sigma$ on $\alpha_n$, and hence on the whole extension $K_\infty$, the character $\chi$ is injective. Since $\Gal(K_\infty/K)$ embeds into an abelian group, it is itself abelian.
\end{proof}
\begin{corollary}\label{CorBergerApp}
	Assume $p$ is odd. Suppose $\alpha_0$ is a uniformizer for $K$, $P'(0)$ is nonzero, and we are given a branch associated to $P(x)$ and $\alpha_0$ which is \textit{potentially} tamely ramification-stable.
	
	If $K_\infty/K$ is Galois, it has a finite-index abelian subgroup.
\end{corollary}
\begin{proof}
	Select $N$ such that the branch is tamely ramification-stable over $K_N$. Since it is still the case that the polynomials $P^n(x) - \alpha_0$ are Eisenstein, the new base point $\alpha_N$ remains a uniformizer. Therefore, Theorem~\ref{ThmBergerApp} applies over this larger field, and hence $\Gal(K_\infty/K_N)$ is abelian. Its index in $\Gal(K_\infty/K)$ is exactly $q^N$.
\end{proof}

We cannot relax the assumption that $\alpha_0$ is a uniformizer, as this is crucial to Berger's construction of the Coleman power series. Moreover, the fact that $\alpha_0$ is a uniformizer means that every $\alpha_n$ will also be a uniformizer of the field it generates over $K$, and so $d=1$ for any branch based at $\alpha_0$. As a result, whether or not the branch is potentially tamely ramification-stable depends \emph{entirely} on $P(x)$.

Theorem~\ref{ThmBergerApp} is not vacuous; there are tamely ramification-stable branches associated to Galois extensions. For example, it is straightforward to check that Berger's example (Theorem~6.5~\cite{BergerIterated}) 
$$P(x) = x^3 + 6x^2 + 9x\quad \alpha_0 = -3\quad K=\Q_3$$
satisfies Theorem~\ref{ThmBergerApp} by combining our observation that $d=1$ with the effective results of Section~\ref{SectEffective}. 

In fact, because $d=1$ and the polynomial in question has prime degree, the branch is guaranteed to be potentially tamely ramification-stable, so we could have applied Corollary~\ref{CorBergerApp}, without making any calculations, to determine that the Galois group has a large abelian subgroup (applying our effective results, one can see that this would have proven $K_\infty/K_1$ is abelian). This can be done for any other examples involving a post-critically bounded or prime degree polynomial.

\subsection{A question about wild ramification in arboreal extensions.}

Both Aitken, Hajir, and Maire~\cite{AitkenHajirMaire} (Question 7.1) and Bridy, Ingram, Jones, Juul, Levy, Manes, Rubinstein-Salzedo, and Silverman~\cite{BIJJLMRS} (Conjecture 6) raise questions about wild ramification in arboreal extensions. Namely: are there arboreal extensions over number fields which are ramified at only finitely many primes but \textit{not} wildly ramified?

We answer this negatively for all arboreal extensions associated to polynomials of prime-power degree. Under some restrictions on the base point, we can also show that such arboreal extensions are not only infinitely wildly ramified, but that all of their higher ramification subgroups are nontrivial. For the latter, we do not need the full strength of our results, only that $K_\infty/K$ is arithmetically profinite (which, for certain base points, already follows from Cais, Davis, and Lubin~\cite{CDL}).

\BIJJLMRSApplication
\begin{proof}
	It suffices to work over the completion $K$ of $F$ at a prime lying over $p$, and we may also take finite extensions of the ground field as necessary. Iteration and conjugation of PCB polynomials are PCB, so we may replace $P(x)$ by some conjugate iterate of itself, which allows us to modify its degree and ensure it fixes $0$. So by Proposition~\ref{PropPCB} we may assume that in addition to being monic, $P(x)$ has integral coefficients, and fixes $0$. Replacing $P$ by $P^s$ for a sufficiently large integer $s$, we may assume that the the size of the residue field of $K$ divides the degree of $P$.

	Recall that our results require $v(\alpha_0) \neq 0$. If $v(\alpha_0) = 0$, then after possibly extending $F$, we will conjugate by a translation to make its valuation positive. In particular, $P(x)$ has a fixed point congruent to $\alpha_0$ modulo $\pi_K$, because
	$$P(x) - x \equiv x^{p^r} - x \mod \pi_K,$$
	and the size of the residue field divides $p^r$, so that every element of the residue field is a zero of $P(x)-x$ modulo $\pi_K$. Let $\alpha$ be such a fixed point, then replace $P(x)$ by its conjugate by $x\mapsto x - \alpha$ and $\alpha_0$ by $\alpha_0 - \alpha$.

	This leaves us with a final pair $(P(x),\alpha_0)$ where $v(\alpha_0) \neq 0$. It follows from Proposition~\ref{PropSemistable} that (every) branch extension $K_\infty/K$ is infinitely wildly ramified, hence the full arboreal extension $K_{arb}/K$ is also infinitely wildly ramified.
	
	Because being post-critically bounded is conjugation and composition invariant, we may \textit{always} assume when $P(x)$ is post-critically bounded that $v_\pp(\alpha_0) \neq 0$.

	We can say more if $P(x)$ has prime degree with $v_\pp(\alpha_0) \neq 0$ or $P(x)$ is post-critically bounded and $v_\pp(\alpha_0) \neq 0$, and there is a branch such that $p\nmid d$, as then Corollary~\ref{CorollaryMain} applies: there is an $N$ such that after replacing $K$ by $K_N$,
	$$K_n = K_\infty^{((V-1)(n-N) + 1)}.$$
	
	Those are the subfields of $K_\infty$ fixed by $\Gamma_K^{b_{(V-1)(n-N) + 1}}$. The branch extension $K_\infty/K$ is contained in the full arboreal extension $K_{arb}/K$, which, combined with the functoriality of the upper numbering, means $K_n$ is the subfield of $K_\infty/K$ that is fixed by the subgroups $\Gamma_{arb}^{b_{(V-1)(n-N) + 1}}$. But the fields $K_n$ are all distinct, and hence the subgroups which fix them must all be distinct too. Finally, it was shown that the ramification breaks $b_{(V-1)(n-N) + 1}$ are unbounded as a function of $n$, and so every upper-numbered higher ramification subgroup of $\Gamma_{arb}$ is nontrivial.
\end{proof}
\begin{observation}
	Bridy, Ingram, Jones, Juul, Levy, Manes, Rubinstein-Salzedo, and Silverman~\cite{BIJJLMRS} showed that a finitely ramified arboreal extension over a number field necessarily comes from a post-critically finite, and hence post-critically bounded map. This means that the preceding theorem applies as soon as one checks that $p$ does not divide $d$ (the stronger case, without restricting $v_\pp(\alpha_0)$ because the map is PCB).
\end{observation}

The theorem tells us that, at least in some cases, the higher ramification subgroups of $\Gamma_{arb}$ are all nontrivial, so we are led to wonder how large or small these subgroups might be. In particular, is $K_{arb}/K$ arithmetically profinite? We suspect not, and conjecture that if there is no branch such that $K_\infty/K$ is Galois, then the wild ramification subgroup has infinite index inside $\Gamma_{arb}$ (in other words, the tame part of $K_\infty/K$ has infinite degree over $K$). However, it seems plausible that this could be the only obstacle to the extension being arithmetically profinite: is it the case that for any $1 < \mu < \nu$, the index $[\Gamma_{arb}^\mu:\Gamma_{arb}^\nu]$ is finite?

\subsection{Effective results; calculating limiting ramification data.}\label{SectEffective}

Almost every step of the proof is effective, and in practice straightforward to compute. Here we sketch the computation of most of the limiting ramification data (Definition~\ref{defnlimitingramdata}). An implementation in SageMath~\cite{Sage} is available upon request. The only ineffective step made to obtain our results occurs in Proposition~\ref{PropSemistable}. The determination of ``sufficiently large'' to ensure that (b) and (c) of this proposition are satisfied is not effective. This also means that the value $d = \lim_{n\rightarrow\infty} d_n$ is not effective. Knowing that $p$ does not divide $d$ is an important input to our main results, so from a computational perspective, this is a particularly unfortunate limitation.

However, if $d$ is known, then all of our constants are effective. For example: if $\alpha_0$ is a uniformizer, such as in the previous section, then $P^n(x)-\alpha_0$ is Eisenstein, so $\alpha_n$ is also a uniformizer, and so $d=1$ and (H) is immediately satisfied at the first level.

\subsubsection{Calculating $V$, $R$, $M$, and $E$.}
We begin with the computation of $V$, $R$, $M$, and $E$: the number of vertices, the (logarithm of) the $x$-coordinates of the vertices, and the main and error terms describing the heights of the vertices. Interestingly, these depend only on the valuations of the coefficients of $P$ and on the \textit{sign} of the valuation of $\alpha_0$. They do not depend on the choice of branch.

All of the following steps can be extracted readily from the proof of Proposition~\ref{PropVertexStability}. Roughly, the proposition tells us that when $v(\alpha_n)$ is small, we can drop the small error terms that show up in the minimum defining the Newton polygon $\NP_n$ as long as we carefully track which terms achieve that minimum.

\begin{enumerate}[]
	\item \textit{Step 1.} For each $0\leq k \leq r$, compute the minimum
	\begin{equation}\label{algmin}
	M_{p^k} = \min_{p^k\leq j\leq q} \left\{v\binom j {p^k} + v(P_j)\right\}.
	\end{equation}
	\item \textit{Step 2.} For each $0\leq k\leq r$: if $v(\alpha_0)$ is positive (resp. negative), let $j$ be the first (resp. last) index achieving the minimum (\ref{algmin}) which defines $M_{p^k}$ . Then set
	$$E_{p^k} = j- p^k.$$
	\item \textit{Step 3.} Let $\NP$ be the lower convex hull of the following vertices:
	\begin{equation*}\label{algNP}
	\left\{\left(p^k, M_{p^k} + E_{p^k} \frac 1 {q^2}\right)\ :\ 0\leq k \leq r\right\}.
\end{equation*}
	The division by $q^2$ is arbitrary - any larger power of $q$ will work as well. This polygon is an approximation to the polygons $\NP_n$ that is precise enough to contain all the limiting ramification data. It is important that the error term is still present, because it can contribute vertices to the Newton polygons even though its contribution decreases rapidly. Degenerating all the way to the convex hull of the points $(p^k, M_{p^k})$ will lose this crucial information.
	
	\item \textit{Step 4.} Let $V$ be the number of vertices of the polygon $\NP$, and write the $x$-coordinates of the vertices of $\NP$ as $p^{r_1},...,p^{r_V}$. Then the limiting ramification data is:
	\begin{align*}
	V(P,\alpha_0) &= V\\
	R(P,\alpha_0) &= (r_1,...,r_V)\\
	M(P,\alpha_0) &= (m_{1},...,m_{V})\\
	E(P,\alpha_0) &= (e_{1},...,e_{V})
\end{align*}
(recall $m_i = M_{p^{r_i}}$, likewise $e_i = E_{p^{r_i}}$).
\end{enumerate}

\subsubsection{Calculating $C$.}
The constant $C$ requires slightly more information to calculate. Unlike $V$, $R$, $M$, and $E$, this constant depends on the branch. However, the dependence is weaker than one might expect: if $\alpha_0 \neq 0$, there is a constant $N$ which is \textit{uniform} among all branches when the valuation of the base point, $v(\alpha_0)$, is fixed, such that $C = q^N v(\alpha_N)$. In fact, this constant $N$ does not even depend on $P(x)$, only its degree. When $\alpha_0 = 0$, there is still such a constant, but it depends on the valuations of the coefficients of $P(x)$ and the number of leading zeros of the branch.

Inspecting the proof of Proposition~\ref{PropVertexStability}, we see that if we have an $N$ such that $(P,\alpha_N)$ satisfy Proposition~\ref{PropSemistable}(a), then the constant $C$ is given by $q^Nv(\alpha_N)$. So we simply need to give a bound on this $N$ in terms of $P$ and $v(\alpha_0)$.

We can extract this from the proof of Proposition~\ref{PropSemistable}(a). If $v(\alpha_0) < 0$ then we are done. If $v(\alpha_0) > 0$, more work is required.

If $\alpha_0 \neq 0$, then the decrease in valuation is partly controlled by the following estimate:
$$v(\alpha_n) \leq \max \{v(\alpha_{n-1})-1, v(\alpha_{n-1})/2\}.$$

In the maximum, it is easy to see that
$$v(\alpha_{n-1})-1\leq v(\alpha_{n-1})/2$$
if and only if
$$v(\alpha_{n-1})\leq 2,$$
and when that occurs, it must be that $v(\alpha_n) \leq 1$. So after $N=v(\alpha_0)$ steps, we are guaranteed to be in a situation where Proposition~\ref{PropSemistable}(a) applies, and hence $C = q^N v(\alpha_N)$.

Otherwise, $\alpha_0 = 0$. Let $k$ be the number of leading $0$s in the branch, which means $\alpha_{k} \neq 0$ and $\alpha_{k-1} = 0$. and by inspecting the Newton polygon of $P(x) - \alpha_{k-1} = P(x)$, a generous bound for $v(\alpha_k)$ is $\ell = \max\{v(P_j)\}$, as long as $\alpha_1 \neq 0$. Then we may apply our reasoning for the case $\alpha_0 \neq 0$, but with $\alpha_k$ in place of $\alpha_0$ to see that
$$C = q^{k+\ell}v(\alpha_{k+\ell})$$

This gives us a remarkably simple process for computing the index $N$ such that $C = q^N v(\alpha_N)$, and of course $C$ itself:
\begin{enumerate}[]
	\item \textit{Step 1.} If $v(\alpha_0) < 0$, then let $N = 0$.
	\item \textit{Step 2.} If $v(\alpha_0) > 0$ and $\alpha_0\neq 0$, then let $N = v(\alpha_0)$.
	\item \textit{Step 3.} If $\alpha_0 = 0$, let $k$ be the number of leading zeros in the branch and let $\ell = \max\{v(P_j)\}$. Then let $N = k+\ell$.
	\item \textit{Step 4.} Set $C = q^N v(\alpha_N)$.
\end{enumerate}

Evidently, the value $N$ is independent of the branch except when $\alpha_0 = 0$, and in that case the dependence is only on the number of leading zeros. Usually this index is much larger than necessary.

\subsubsection{Sample calculation.}

In any particular case, it is almost straightforward to check that a pair is tamely ramification-stable, \textit{except} for the tameness component, since we do not have an effective way to compute $d$. However, it is still possible to do so in some cases.

The following example is small enough that one can carry out the calculation by hand.

Let $K=E=\Q_3(\sqrt 3)$ with valuation $v$ normalized so that $v(\sqrt 3) = 1$. Consider the polynomial
\begin{align*}
	P(x) &= x^9 + 12\sqrt 3x^7 + 18x^6 + 3\sqrt 3x^4 + \frac 3 5x^3 + 9x,
\end{align*}
with any branch whose initial sequence of valuations looks like $(4,2/3,2/27,...)$.

We readily obtain our effective constants:
\begin{align*}
	V &= 3,\\
	R &= (0,1,2),\\
	M &= (3,2,0),\\
	E &= (3,0,0).
\end{align*}

as well as
$$C = 9^4 v(\alpha_4) = 9^4 * \frac 2 3 \frac 1 {9^3} = 6.$$

Inspecting the first few levels of such a branch in Sage, we see that our sequence $d_n$ looks like $4,2,2,...$, hence $d=2$, which is not divisible by $p=3$. To be more precise, while the value $N$ from~\ref{PropSemistable} is not effectively determined, we can see from the proof that as soon as some $d_n$ is not divisible by $p$ and $v(\alpha_n) < 1$, we have reached a suitable index. This is because at each step there is no way for the valuation to decrease by a factor of $q$ without the ramification index being $q$ as well. Combined with this limiting ramification data, one can see directly that $(P,\alpha_1)$ is tamely ramification-stable. Therefore, when we replace $K$ by $K_1$, we may apply Theorem~\ref{TheoremMain} to obtain
$$K_n = K_\infty^{((V-1)(n-1) + 1)}.$$

So, even though $P(x)$ is not prime-degree or post-critically bounded, it is an example of a polynomial whose branch extensions are amenable to study by our methods.

\section*{Acknowledgments}
I would like to thank my advisor, Joseph H. Silverman, for many helpful discussions on this project and his careful comments on early versions of the paper. I would also like to thank the anonymous referee for their comments, which greatly improved the clarity of the paper.

\end{document}